\PassOptionsToPackage{quiet}{fontspec}
\documentclass[leqno]{article}

\usepackage{amsmath,amsthm,amsfonts,color,upgreek}
\usepackage{amssymb}
\usepackage{graphicx}
\usepackage{setspace}
\usepackage{mathrsfs}
\usepackage{graphicx}
\usepackage[shortlabels]{enumitem}
\usepackage{authblk}
\usepackage{verbatim}

\newtheorem{theorem}{Theorem}

\newtheorem{lemma}{Lemma}
\newtheorem{proposition}{Proposition}
\newtheorem{example}{Example}

\newtheorem{case}{Case}

\newtheorem{claim}{Claim}
\newtheorem{corollary}{Corollary}
\newtheorem{definition}{Definition}
\newtheorem{problem}{Problem}
\newtheorem{remark}{Remark}

\usepackage{chngcntr}
\numberwithin{equation}{section}

\usepackage{calc}
\usepackage{tikz}
\usepackage{geometry}
\usepackage{pgffor}
\usepackage{ifthen}
\usetikzlibrary{arrows.meta}
\usetikzlibrary{decorations.markings}
\tikzstyle{vertex}=[circle, draw, inner sep=2pt, minimum size=6pt]
\tikzstyle{filledvertex}=[circle, draw, fill, inner sep=2pt, minimum size=6pt]

\let\leq\leqslant
\let\geq\geqslant

\begin{document}
\begin{spacing}{1.15}

\title{Cross-intersecting families with covering number constraints}
\author{Yandong Bai
\thanks{School of Mathematics and Statistics, 
Xi’an-Budapest Joint Research Center for Combinatorics,
Northwestern Polytechnical University,
Xi'an, Shaanxi 710129, China;
Research $\&$ Development Institute of 
Northwestern Polytechnical University in Shenzhen,
Shenzhen, Guangdong 518057, China.
E-mail: {\tt bai@nwpu.edu.cn}.}
\quad Haoyun Gu
\thanks{School of Mathematics and Statistics, 
Xi’an-Budapest Joint Research Center for Combinatorics,
Northwestern Polytechnical University,
Xi'an, Shaanxi 710129, China.
E-mail: {\tt haoyungu@nwpu.edu.cn}.}
}
\date{\today}

\maketitle

\begin{abstract}
    Two families $\mathcal{F}$ and $\mathcal{G}$ are cross-intersecting if every set in $\mathcal{F}$ intersects every set in $\mathcal{G}$. The covering number $\tau(\mathcal{F})$ of a family $\mathcal{F}$ is the minimum size of a set that intersects every member of $\mathcal{F}$. In 1992, Frankl and Tokushige determined the maximum of $|\mathcal{F}| + |\mathcal{G}|$ for cross-intersecting families $\mathcal{F} \subset \binom{[n]}{a}$ and $\mathcal{G} \subset \binom{[n]}{b}$ that are non-empty (covering number at least 1) and also characterized the extremal configurations. This seminar result was recently extended by Frankl (2024) and Frankl and Wang (2025) to cases where both families are non-trivial (covering number at least 2), and where one is non-empty and the other non-trivial, respectively. In this paper, we establish a unified stability hierarchy for cross-intersecting families under general covering number constraints. We determine the maximum of $|\mathcal{F}| + |\mathcal{G}|$ for cross-intersecting families $\mathcal{F} \subset \binom{[n]}{a}$ and $\mathcal{G} \subset \binom{[n]}{b}$ with the following covering number constraints: 
        (1) $\tau(\mathcal{F}) \geq s$ and $\tau(\mathcal{G}) \geq t$; 
        (2) $\tau(\mathcal{F}) = s$ and $\tau(\mathcal{G}) \geq t \geq 2$; 
        (3) $\tau(\mathcal{F}) \geq s$ and $\tau(\mathcal{G}) = t$;
        (4) $\tau(\mathcal{F}) = s$ and $\tau(\mathcal{G}) = t$;
    provided $a \geq b + t - 1$ and $n \geq \max\{a + b, bt\}$. The corresponding extremal families achieving the upper bounds are also characterized.
\end{abstract}

\section{Introduction}

    Let $[n] := \{1, \ldots, n\}$ denote the standard $n$-element set and $[p,q]$ the set $\{p,p+1,\ldots,q\}$ for $q\geq p$. We write $2^{[n]}$ for the power set of $[n]$ and $\binom{[n]}{k}$ for the collection of all $k$-subsets of $[n]$. A family $\mathcal{F}\subset 2^{[n]}$ is \emph{intersecting} if any two of its members have non-empty intersection, and it is {\em $k$-uniform} if $\mathcal{F}\subset \binom{[n]}{k}$. Two families $\mathcal{F}$ and $\mathcal{G}$ are \emph{cross-intersecting} if $F \cap G \neq \emptyset$ for all $F \in \mathcal{F}$ and $G \in \mathcal{G}$. 
    
    By definition, every family and itself are cross-intersecting. So the concept cross-intersecting can be viewed as a generalization of intersecting. Besides, in many extremal problems concerning intersecting families, a frequently employed technique is to decompose the family into subfamilies and then, by analyzing some cross-intersecting subfamilies, to bound the size of the original family. Therefore, the study of cross-intersecting families has attracted considerable attention in extremal set theory.
    
    A fundamental parameter associated with a family $\mathcal{F} \subset 2^{[n]}$ containing no empty set is its \emph{covering number} $\tau(\mathcal{F})$, defined as the minimum size of a set $T \subset [n]$ that intersects every member of $\mathcal{F}$, i.e., 
    \[
       \tau(\mathcal{F})
       :=\min\{|T|: T\subset [n], T\cap F\neq \emptyset \text{~for all~} F\in \mathcal{F}\}.
    \]
    Note that a family $\mathcal{F}$ containing no empty set is non-empty precisely when $\tau(\mathcal{F}) \geq 1$, and is \emph{non-trivial} (i.e., has empty total intersection) exactly when $\tau(\mathcal{F}) \geq 2$. From this point of view, the covering number naturally induces a hierarchy for classifying intersecting families.
    
    In what follows, we first briefly review key results on maximum sizes of intersecting uniform families under covering number constraints, then turn to cross-intersecting families. By establishing a stability hierarchy, we extend the aforementioned results and provide a comprehensive framework for cross-intersecting uniform families with general covering number constraints.

\subsection{Intersecting families with covering number constraints}
    
    Extremal set theory is concerned with determining the maximum possible size of a set system satisfying certain constraints. A foundational result in this area is the Erdős–Ko–Rado (EKR) theorem, which establishes such a bound for $k$-uniform intersecting families. 
    
    \begin{theorem}[Erdős, Ko and Rado \cite{EKR61}]\label{thm: EKR}
        Suppose that $\mathcal{F} \subset \binom{[n]}{k}$ is an intersecting family with $n \geq 2k\geqslant 4$. Then
        \[
            |\mathcal{F}| \leq \binom{n-1}{k-1}. 
        \]
    \end{theorem}
    
    For $n > 2k$, the equality holds if and only if $\mathcal{F}$ is a \emph{full star}, i.e., all $k$-sets containing a fixed element,
    \[
        \mathcal{S}_i(n,k)
        :=\left\{F\in \binom{[n]}{k}: i\in [n]\right\}.
    \]

    Since the empty set is not contained in any uniform family, the EKR theorem characterizes the largest intersecting $k$-uniform families with covering number at least 1. Stability results for the EKR theorem explore the structure and maximum possible size of intersecting families that are different from the extremal configuration, namely, the full star. A prominent line of research concerns classifying these families based on their covering numbers. Hilton and Milner \cite{HM67} initiated this direction by demonstrating the largest $k$-uniform intersecting families with covering number at least 2. 
    
    \begin{theorem}[Hilton and Milner \cite{HM67}]\label{thm: HM}
        Suppose that $\mathcal{F} \subset \binom{[n]}{k}$ is an intersecting family with $n > 2k\geqslant 4$ and $\tau(\mathcal{F})\geq 2$. Then
        \[
            |\mathcal{F}| \leq \binom{n-1}{k-1} - \binom{n-k-1}{k-1} + 1. 
        \]
        Moreover, for $n > 2k$, the equality holds if and only if $\mathcal{F}$ is isomorphic to the following family,
        \[
        \mathcal{H}(n,k)
        :=\left\{F\in \binom{[n]}{k}:1\in F,~F\cap [2,k+1]\neq \emptyset\right\}\cup \{[2,k+1]\}.
        \]
        for $k=3$, there is one more possibility, namely, the triangle family,
        \[
        \mathcal{K}(n,k)
        :=\left\{F\in \binom{[n]}{k}:1\in F, |F\cap [3]|\geqslant 2\right\}.
        \]
    \end{theorem}

Let $\mathcal{G}(n,k):=\mathcal{A}\cup \mathcal{B}$, where
\[
\mathcal{A} 
:= \left\{A \in \binom{[n]}{k} : 1 \in A \text{ and } A \cap B \neq \emptyset \text{ for all } B \in \mathcal{B} \right\},~
\mathcal{B} 
:= \{[2,k+1], \{2\} \cup [k+2,2k], \{3\} \cup [k+2,2k]\}.
\]
    
    For intersecting families with covering number at least 3, the construction $\mathcal{G}(n,k)$ gives a lower bound, and it was proved by Frankl \cite{F80} that for $k \geq 4$ and $n > n_0(k)$ (with $n_0(k)$ exponential in $k$), this family is optimal and up to isomorphism, $\mathcal{G}(n,k)$ is the unique optimal family. Recently, Frankl and Wang~\cite{FW25} significantly improve the range of $n$ for which the extremal result holds. 
   
\begin{theorem}[Frankl and Wang~\cite{FW25}]\label{thm: covering number 3}
Suppose that $\mathcal{F} \subset \binom{[n]}{k}$ is an intersecting family with $n \geq 2k\geqslant 14$ and $\tau(\mathcal{F})\geq 3$. Then
\[
|\mathcal{F}| \leq \binom{n-1}{k-1} - \binom{n-k}{k-1} - \binom{n-k-1}{k-1} + \binom{n-2k}{k-1} + \binom{n-k-2}{k-3} + 3.
\]
Moreover, for $n > 2k$, the equality holds if and only if $\mathcal{F}$ is isomorphic to $\mathcal{G}(n,k)$.
\end{theorem}  

We mention that Kupavskii \cite{K24} independently obtain the same statements for $k\geq 100$ and $n>2k$ by a different proof. 

Define the function
\begin{equation}\label{equation: f(n,k,s)}
    m(n,k,s) := \max \left\{ |\mathcal{F}| : \mathcal{F} \subset \binom{[n]}{k} \text{ is intersecting and } \tau(\mathcal{F}) \geq s \right\}.
\end{equation}
With this terminology, Theorems \ref{thm: EKR}, \ref{thm: HM} and \ref{thm: covering number 3} can be stated as
    \[
       m(n,k,1) = \binom{n-1}{k-1} \text{ for } n \geq 2k\geqslant 4;
    \]
    \[
       m(n,k,2) = \binom{n-1}{k-1} - \binom{n-k-1}{k-1} + 1 \text{ for } n > 2k\geqslant 4;
    \]
    \[
       m(n,k,3) = \binom{n-1}{k-1} - \binom{n-k}{k-1} - \binom{n-k-1}{k-1} + \binom{n-2k}{k-1} + \binom{n-k-2}{k-3} + 3 \text{ for } n \geqslant 2k\geqslant 14.
    \]
The study of maximum intersecting families with larger covering numbers has also seen many progresses. In particular, $m(n,k,4)$ was determined  by Frankl, Ota and Tokushige \cite{FOT96} for $k \geq 9$ and sufficiently large $n$, by Chiba et al. \cite{C+12} for $k=4$, and by Furuya and Takatou \cite{FT13} for $k=5$. Recently, Frankl and Wang \cite{FWic25} determined $m(n,k,5)$ for $k \geq 69$ and $n \geq 5k^6$. Note that $s\leqslant k$. For the special case $s=k$, Erdős and Lovász \cite{EL75} established the bounds 
\[
\lfloor k!(e-1) \rfloor \leq m(n,k,k) \leq k^k.
\] 
and determined that $f(n,3,3) = 10$.  
Lovász \cite{L75} conjectured that $m(n,k,k)=\lfloor k!(e-1) \rfloor$. This conjecture was disproved for $k \geq 4$ by Frankl, Ota, and Tokushige \cite{FOT96}, who improved the lower bound to
\[
m(n,k,k) \geq 
\begin{cases}
\left(\dfrac{k}{2} + 1\right)^{k-1} & \text{for even } k, \\
\left(\dfrac{k+3}{2}\right)^{(k-1)/2}\left(\dfrac{k+1}{2}\right)^{(k-1)/2} & \text{for odd } k.
\end{cases}
\]
But there is not too many progresses on $m(n,k,s)$ for $5< s<k-1$. Besides, stability problems for the EKR theorem and the Hilton--Milner theorem under other additional constraints have also been extensively studied, such as the maximum degree, the diversity of intersecting families or restrictions only excluding subfamilies of the aforementioned extremal configurations. We refer the reader to \cite{Frankl25, GXZ24, HK17, KM17, K25} for more details. 
        
\subsection{Cross-intersecting families with covering number constraints}
      
    It is natural to introduce the covering number constraints into the study of cross-intersecting families. In 1992, Frankl and Tokushige \cite{FT92} determinied the maximum total size for non-empty cross-intersecting families (i.e., with covering numbers at least 1).

\begin{theorem}[Frankl and Tokushige \cite{FT92}] \label{thm: (1,1) CI}
    Suppose that $\mathcal{F}\subset \binom{[n]}{a}$, $\mathcal{G}\subset \binom{[n]}{b}$ are cross-intersecting families with $a\geq b\geq 2$ and $\tau(\mathcal{F}),\tau(\mathcal{G})\geq 1$. Then for $n \geq a+b$,
    \begin{equation*} 
       |\mathcal{F}| + |\mathcal{G}| 
       \leq \binom{n}{a} - \binom{n-b}{a} + 1.
    \end{equation*} 
    Moreover, unless $n=a+b$ or $a=b=2$, the equality holds if and only if $\mathcal{F}\cup \mathcal{G}$ is isomorphic to $\mathcal{M}^1(n,a,b)$ (see definition \ref{def: M(a,b,t)}).
\end{theorem}


\begin{definition}\label{def: M(a,b,t)}
    Let $\mathcal{B} :=\{B_1,...,B_t \}\subset \binom{[n]}{b}$ such that $B_i\cap B_j = \emptyset$ for all $i\ne j$. Note that $\tau(\mathcal{B})=t$. Denote by
    \begin{equation}
       \mathcal{M}^t(n,a,b) 
       :=\{B_1,,...,B_t \} \cup 
       \left\{A\in \binom{[n]}{a}: A\cap B_i \ne \emptyset \text{ for all } 1\leq i\leq t\right\}.
    \end{equation} 
\end{definition}

In 2024, Frankl \cite{F24} extended Theorem \ref{thm: (1,1) CI} to families that have covering numbers at least 2.
    
\begin{theorem}[Frankl \cite{F24}] \label{thm: (2,2) CI}
    Suppose that $\mathcal{F}\subset \binom{[n]}{a}$, $\mathcal{G}\subset \binom{[n]}{b}$ are cross-intersecting families with $a\geq b\geq 2$ and $\tau(\mathcal{F}),\tau(\mathcal{G})\geq 2$. Then for $n>a+b$,
    \begin{equation*} 
      |\mathcal{F}| + |\mathcal{G}| 
      \leq \binom{n}{a} - 2\binom{n-b}{a} + \binom{n-2b}{a} + 2.
    \end{equation*} 
    Moreover, unless $a=b=2$, the equality holds if and only if $\mathcal{F}\cup \mathcal{G}$ is isomorphic to $\mathcal{M}^2(n,a,b)$.
\end{theorem}

As a key step in proving Theorem \ref{thm: covering number 3}, Frankl and Wang \cite{FW25} further investigated the problem for families such that one has covering number at least 1 and the other at least 2.

\begin{theorem}[Frankl and Wang \cite{FW25}]\label{thm: (1,2) CI}
    Suppose that $\mathcal{F} \subset \binom{[n]}{a}$, $\mathcal{G} \subset \binom{[n]}{b}$ are cross-intersecting families with $b\geq 1$, $a\geq b+1$, $\tau(\mathcal{F})\geqslant 1$ and $\tau(\mathcal{G})\geqslant 2$. Then for $n \geq a+b$, 
    \[
    |\mathcal{F}| + |\mathcal{G}| \leq \binom{n}{a} - 2\binom{n-b}{a} + \binom{n-2b}{a} + 2,
    \]
    Moreover, when $n > a+b$ and $a \geq 2$, the equality holds if and only if $\mathcal{F} \cup \mathcal{G}$ is isomorphic to $\mathcal{M}^2(n,a,b)$.
\end{theorem}

We remark that for cross-intersecting families $\mathcal{F} \subset \binom{[n]}{a}$, $\mathcal{G} \subset \binom{[n]}{b}$ with $a\geq b$ and $n\geq a+b$ but without covering number constraints, the sum of their sizes have a trivial bound 
$
|\mathcal{F}| + |\mathcal{G}| \leq \binom{n}{a}.
$
When $n>a+b$ or $a>b$, the equality holds if and only if $\mathcal{F}=\binom{[n]}{a}$ and $\mathcal{G} = \emptyset$; when $a=b=n/2$, equality can also hold if  $\mathcal{F} = \mathcal{G}$ is a full star $\mathcal{S}_i(n,a)$ or a family such that $|\mathcal{F}\cap \{A,[n]\backslash A\}|=1$ for any $A\subset [n]$ with $|A|=n/2$, this is also the extremal family of the EKR theorem for $n=2k$. 

\subsection{Our main results}

Our work establishes a stability hierarchy for cross-intersecting families under general covering number constraints, unifying and extending previous results in several directions:
\begin{itemize}
    \item We unify and extend Theorems \ref{thm: (1,1) CI} and \ref{thm: (1,2) CI} by determining the maximum sum $|\mathcal{F}| + |\mathcal{G}|$ for cross-intersecting families under the condition $\tau(\mathcal{F}) \geq 1$ and $\tau(\mathcal{G}) \geq t$ for all $t \geq 2$ (Theorem \ref{thm: (1,t) CI}).
    \item We further generalize this to the case where $\tau(\mathcal{F}) = s$, $\tau(\mathcal{G}) \geq t$ or $\tau(\mathcal{F}) = s$, $\tau(\mathcal{G}) = t$ (Theorem 8), and derive a corollary for the cases $\tau(\mathcal{F}) \geq s$, $\tau(\mathcal{G}) \geq t$, or $\tau(\mathcal{F}) \geq s$, $\tau(\mathcal{G}) = t$ (Corollary \ref{corollary: geq (s,t) CI1}).
    \item We determine the maximum sum for initial (left-shifted) cross-intersecting families (Theorem \ref{thm: initial CI}), which will serve as an important tool in the proof of Theorem \ref{thm: (1,t) CI}. Moreover, we obtain the maximum sum for initial cross-intersecting families under the following conditions: (1) $\tau(\mathcal{F}) = s$, $\tau(\mathcal{G}) \geq t \geq 2$; (2) $\tau(\mathcal{F}) = s$, $\tau(\mathcal{G}) = t$; (3) $\tau(\mathcal{F}) \geq s$, $\tau(\mathcal{G}) \geq t$; and (4) $\tau(\mathcal{F}) \geq s$, $\tau(\mathcal{G}) = t$ (Corollaries \ref{corollary: initial CI} and \ref{corollary: initial (s,t) CI}).
\end{itemize}

Let us introduce our results precisely in the following.

\begin{theorem}\label{thm: (1,t) CI} 
    Suppose that $\mathcal{F}\subset \binom{[n]}{a}$, $\mathcal{G}\subset \binom{[n]}{b}$ are cross-intersecting families with $a\geq b + t - 1$, $b\geq 2$, $\tau(\mathcal{F})\geq 1$ and $\tau(\mathcal{G})\geq t$. Then for $n \geq \max \{a+b, bt\}$,
    \begin{equation} \label{mneeq}
    |\mathcal{F}| + |\mathcal{G}| \leq \binom{n}{a} + \sum_{i=1}^t(-1)^i\binom{t}{i}\binom{n-ib}{a} + t.
    \end{equation}
    Moreover, when $n>a+b$, the equality holds if and only if $\mathcal{F}\cup \mathcal{G}$ is isomorphic to $\mathcal{M}^t(n,a,b)$.
\end{theorem}

Note that Theorems \ref{thm: (1,1) CI} and \ref{thm: (1,2) CI} correspond to the cases $t=1$ and $t=2$ in Theorem \ref{thm: (1,t) CI}, respectively.
We further extend Theorem \ref{thm: (1,t) CI} to more general covering number constraints.

\begin{theorem}\label{thm: (s,t) CI} 
    Suppose that $\mathcal{F}\subset \binom{[n]}{a}$, $\mathcal{G}\subset \binom{[n]}{b}$ are cross-intersecting families with $a\geq b + t - 1$ and $b\geq 2$. If (1) $\tau(\mathcal{F})=s$ and $\tau(\mathcal{G})\geq t\geq 2$, or (2) $\tau(\mathcal{F})=s$ and $\tau(\mathcal{G})= t\geq 2$, then for $n \geq \max \{a+b,bt\}$,
    \begin{equation}
    |\mathcal{F}| + |\mathcal{G}| 
    \leq \sum_{i=1}^s \left(\binom{n-i}{a-1} + \sum_{j=1}^{t-1}(-1)^j \binom{t-1}{j} \binom{n-jb-i}{a-1}\right) + \binom{n-s}{b-s} + t-1.
    \end{equation}
    Moreover, when $n>a+b$, the equality holds if and only if $\mathcal{F}\cup \mathcal{G}$ is isomorphic to $\mathcal{M}_s^t(n,a,b)$ (see Definition \ref{def: Ms(a,b,t)}).
\end{theorem}

\begin{definition}\label{def: Ms(a,b,t)}
Let $\mathcal{B}_0 := \{ B \in \binom{[n]}{b}: [1,s] \subset B \} $ and $\mathcal{B}_1 := \{D_1,...,D_{t-1} \}$ such that $D_i \cap [1,s] = \emptyset$ and $D_i \cap D_j = \emptyset, i\ne j$ for all $i,j \in [1,t-1]$. Define $\mathcal{B}' := \mathcal{B}_0 \cup \mathcal{B}_1$. Denote by
\begin{equation}
    \mathcal{M}_s^t(n,a,b) 
    := \left\{ A \in \binom{[n]}{a}: A \cap B \ne \emptyset \text{ for all } B \in \mathcal{B}' \right\} \cup \mathcal{B}'.
\end{equation}   
\end{definition}

Note that $\mathcal{M}_b^t(n,a,b) \cong \mathcal{M}^t(n,a,b)$.
As a corollary, we obtain the following result.

\begin{corollary}\label{corollary: geq (s,t) CI1}
    Suppose that $\mathcal{F}\subset \binom{[n]}{a}$, $\mathcal{G}\subset \binom{[n]}{b}$ are cross-intersecting families with $a\geq b + t - 1$ and $b\geq 2$. 
    If (1) $\tau(\mathcal{F}) \geq s\geq 2$ and $\tau(\mathcal{G}) \geq t\geq 2$; or (2) $\tau(\mathcal{F}) \geq s\geq 2$ and $\tau(\mathcal{G}) = t\geq 2$, then for $n \geq \max \{a+b,bt\}$,
    \begin{equation}
    |\mathcal{F}| + |\mathcal{G}| 
    \leq \binom{n}{a} + \sum_{i=1}^t(-1)^i\binom{t}{i}\binom{n-ib}{a} + t.
    \end{equation}
    Moreover, when $n> a+b$, the equality holds if and only if $\mathcal{F}\cup \mathcal{G}$ is isomorphic to $\mathcal{M}^t(n,a,b)$.
\end{corollary}

\begin{remark}\label{remark: a>b+t-2}
The assumption $a \geq b+t-1$ in Theorems \ref{thm: (1,t) CI} and \ref{thm: (s,t) CI}, and Corollary \ref{corollary: geq (s,t) CI1}  cannot be dropped. If $t\leq a < b+t-1$, then $\mathcal{M}^t(n,a,b)$ may not achieve the maximum of $|\mathcal{F}|+|\mathcal{G}|$. We give some examples. Let $B := \{ b_1,\ldots,b_t \}\subset [n]$. Denote by
\[
  \mathcal{B}' := \left\{ G \in \binom{[n]}{b}: G\cap B\neq \emptyset \right\},
\]
\[
    \mathcal{A}' := \left\{ F \in \binom{[n]}{a}: F \cap G \ne \emptyset \text{ for all } G \in \mathcal{B}' \right\}.
\]
We have $|\mathcal{A}' \cup \mathcal{B}'| = c_1 n^{b-1} + o(n^{b-1})$ and $|\mathcal{M}^t(n,a,b)| = c_2 n^{a-t} + o(n^{a-t})$ for some constant $c_1,c_2$. 
By $a \leq b+t-2$, we get $a-t < b-1$. When $n$ is sufficiently large, it is not difficult to see that $|\mathcal{A}' \cup \mathcal{B}'| > |\mathcal{M}^t(n,a,b)|$.
\end{remark}

\begin{proposition}
Suppose that $\mathcal{F}\subset \binom{[n]}{a}$, $\mathcal{G}\subset \binom{[n]}{b}$ are cross-intersecting families with $t \leq a \leq b + t - 2$, $b\geq 2$, $\tau(\mathcal{F})\geq 1$ and $\tau(\mathcal{G})\geq t\geq 2$. Then for sufficiently large $n$,
    \begin{equation} \label{1tseq}
    |\mathcal{F}| + |\mathcal{G}| \leq \binom{n}{b} - \binom{n-a}{b} + 1.
    \end{equation}    
\end{proposition}

\begin{proof}
By induction on $t$. 
If $t = 1$, then (\ref{1tseq}) follows from Theorem \ref{thm: (1,1) CI}. 
Assume that $t\geq 2$ and (\ref{1tseq}) holds for all $t' < t$.
Note that $|\mathcal{F}| + |\mathcal{G}| = |\mathcal{F}(1)| + |\mathcal{G}(\overline{1})| + |\mathcal{F}(\overline{1})| + |\mathcal{G}(1)|$. 
 
Since $\mathcal{F}(1),\mathcal{G}(\overline{1})$ are cross-intersecting families with $\tau(\mathcal{G}(\overline{1})) \geq t-1$ and $a-1 \leq b + (t-1)-2$, by induction hypothesis, we have 
$|\mathcal{F}(1)| + |\mathcal{G}(\overline{1})| \leq \binom{n-1}{b} - \binom{n-a}{b} + 1.
$

Assume first that $\mathcal{F}(\overline{1}) \ne \emptyset$. Then 
$$
\mathcal{G}(1) = \{ G \in \binom{[2,n]}{b-1}: F \subset G, F \in \mathcal{T}(\mathcal{F}(\overline{1})) \},~~ 
\mathcal{F}(\overline{1}) = \{ F \in \binom{[2,n]}{a}: F \subset G, G \in \mathcal{T}(\mathcal{G}(1) \cup \mathcal{G}(\overline{1})) \}. 
$$
Since $\tau(\mathcal{G}(1) \cup \mathcal{G}(\overline{1})) \geq t$, then we have $|\mathcal{G}(1) \cup \mathcal{G}(\overline{1})| \geq t$ and we can find $t$ different sets $G_1,...,G_t$ in $\mathcal{G}(1) \cup \mathcal{G}(\overline{1})$. For a covering set $S = \{s_1,...,s_t \} \in \mathcal{T}^{(t)}(\mathcal{G}(1) \cup \mathcal{G}(\overline{1}))$, each $s_i$ has at most $|G_i| \leq b$ choices, implying $|\mathcal{T}^{(t)}(\mathcal{G}(1) \cup \mathcal{G}(\overline{1}))| \leq b^t$ (the same for $|\mathcal{T}^{(1)}(\mathcal{F}(\overline{1}))|$). 

Now we can get that $|\mathcal{G}(1)| + |\mathcal{F}(\overline{1})| \leq c \binom{n-1}{b-2} + o(n^{b-2})$, $c \leq b^t+a$. When $n$ is sufficiently large, we have $|\mathcal{G}(1)| + |\mathcal{F}(\overline{1})| < \binom{n-1}{b-1}$, implying $|\mathcal{F}| + |\mathcal{G}| < \binom{n-1}{b} - \binom{n-a}{b} + 1 + \binom{n-1}{b-1} = \binom{n}{b} - \binom{n-a}{b} + 1$, as desired. 
Now assume that $\mathcal{F}(\overline{1}) = \emptyset$. Since $|\mathcal{G}(1)| \leq \binom{n-1}{b-1}$, we have $|\mathcal{F}| + |\mathcal{G}| \leq \binom{n}{b} - \binom{n-a}{b} + 1$, as desired.
\end{proof}

Notice that (\ref{1tseq}) may not hold for small $n$. When $n=40$, $a=t=10$, $b=3$, we have $|\mathcal{M}^{10}(40,10,3)| > 3^{10} = 59049$, while one can see that $\binom{n}{b} - \binom{n-a}{b} + 1 < \binom{n}{b} = \binom{40}{3} = 2080 < 59049$.


\paragraph*{Shifting operation}

    Let $ 1 \leq i < j \leq n $. The \emph{shifting operation} $ S_{ij} $ on $\mathcal{F}$ is defined as follows, 
    \[
    S_{ij}(\mathcal{F}) := \left\{ S_{ij}(F) : F \in \mathcal{F} \right\}, \text{ where}
    \]
    \[
    S_{ij}(F) := 
    \begin{cases}
    (F \setminus \{j\}) \cup \{i\}, & \text{if } j \in F,~i \notin F \text{ and } (F \setminus \{j\}) \cup \{i\} \notin \mathcal{F}; \\
    F, & \text{otherwise}.
    \end{cases}
    \]
This operation replaces $ j $ with $ i $ in $ F $ whenever possible, without creating duplicates. Shifting preserves many properties, such as the intersecting or cross-intersecting property, but may decrease the covering number $ \tau(\mathcal{F})$. This is also one of the main challenges we need to overcome when we use shifting operations in our proofs. 
    A family is \emph{initial} (or \emph{left-shifted}) if $S_{ij}(\mathcal{F})=\mathcal{F}$ for all $ 1 \leq i < j \leq n $.

Define
\[
    \mathcal{H}^t(n,a,b) 
    :=\left\{ A \in \binom{[n]}{a}: |A \cap [b+t-1]| \geq t \right\} 
    \cup \binom{[b+t-1]}{b}.
\]
Note that
\begin{equation} \label{icime}
    |\mathcal{H}^t(n,a,b)|
    =\binom{n}{a} - \sum_{i=0}^{t-1} \binom{b+t-1}{i} \binom{n-b-t+1}{a-i} + \binom{b+t-1}{b}.
\end{equation}

Define
\[
    \mathcal{H}_s^t(n,a,b) 
    :=\left\{ A \in \binom{[n]}{a}: A \cap [1,s] \ne \emptyset, |A \cap [b+t-1]| \geq t \right\} \cup \left\{ B \in \binom{[n]}{b}: [1,s] \subset B \right\}
    \cup \binom{[b+t-1]}{b}.
\]
Note that $\mathcal{H}^t(n,a,b) \cong \mathcal{H}_b^t(n,a,b)$ and 
\begin{equation} \label{icime st}
   |\mathcal{H}_s^t(n,a,b)|
   =\sum_{j=1}^s \left( \binom{n-j}{a-1} - \sum_{i=0}^{t-j-1} \binom{b+t-j-1}{i} \binom{n-b-t+1}{a-i-1} \right) + \binom{n-s}{b-s} + \binom{b+t-1}{b} - \binom{b+t-s-1}{b-s}.
\end{equation}
    
We determine the maximum sum of the sizes of two initial cross-intersecting uniform families. This will also be used in the proofs of Theorems \ref{thm: (1,t) CI} and \ref{thm: (s,t) CI}.

\begin{theorem} \label{thm: initial CI}
    Suppose that $\mathcal{F} \subset \binom{[n]}{a}$, $\mathcal{G} \subset \binom{[n]}{b}$ are initial and cross-intersecting families with $a \geq b+t-1 $ and $\tau(\mathcal{F})\geq 1$. Then for $n\geq a+b$,
    \begin{equation*} 
        |\mathcal{F}| + |\mathcal{G}| 
        \leq |\mathcal{H}^t(n,a,b)|.
    \end{equation*}
    Moreover, when $n > a+b$, the equality holds if and only if $\mathcal{F} \cup \mathcal{G}$ is isomorphic to $\mathcal{H}^t(n,a,b)$.
\end{theorem}

\begin{theorem} \label{thm: initial st CI}
Suppose that $\mathcal{F} \subset \binom{[n]}{a}$, $\mathcal{G} \subset \binom{[n]}{b}$ are initial and cross-intersecting families with $a \geq b+t-1 $, $\tau(\mathcal{F}) = s$ and $\tau(\mathcal{G}) \geq t \geq 2$. Then for $n\geq a+b$,
\begin{equation*} 
    |\mathcal{F}| + |\mathcal{G}|\leq |\mathcal{H}_s^t(n,a,b)|.
\end{equation*}
Moreover, when $n > a+b$, the equality holds if and only if $\mathcal{F} \cup \mathcal{G}$ is isomorphic to $\mathcal{H}_s^t(n,a,b)$.
\end{theorem}

Consequently, we obtain the following results.

\begin{corollary}\label{corollary: initial CI}
    Suppose that $\mathcal{F}\subset \binom{[n]}{a}$, $\mathcal{G}\subset \binom{[n]}{b}$ are initial cross-intersecting families with $a\geq b + t - 1$ and $b\geq 2$. If (1) $\tau(\mathcal{F}) \geq s\geq 2$ and $\tau(\mathcal{G}) \geq t\geq 2$, or (2) $\tau(\mathcal{F})\geq s \geq 2$ and $ \tau(\mathcal{G})= t\geq 2$, then for $n \geq a+b$,
    \begin{align*}
    |\mathcal{F}| + |\mathcal{G}| 
    \leq |\mathcal{H}^t(n,a,b)|.
    \end{align*}
    Moreover, when $n>a+b$, the equality holds if and only if $\mathcal{F}\cup \mathcal{G}$ is isomorphic to $\mathcal{H}^t(n,a,b)$.
\end{corollary}

\begin{corollary}\label{corollary: initial (s,t) CI}
    Suppose that $\mathcal{F}\subset \binom{[n]}{a}$, $\mathcal{G}\subset \binom{[n]}{b}$ are initial cross-intersecting families with $a\geq b + t - 1$ and $b\geq 2$. If (1) $\tau(\mathcal{F})= s\geq 2$ and $\tau(\mathcal{G})\geq t\geq 2$, or (2) $\tau(\mathcal{F})= s\geq 2$ and $\tau(\mathcal{G})= t\geq 2$, then for $n \geq a+b$,
    \begin{align*}
    |\mathcal{F}| + |\mathcal{G}| 
    \leq |\mathcal{H}_s^t(n,a,b)|.
    \end{align*}
    Moreover, when $n>a+b$, the equality holds if and only if $\mathcal{F}\cup \mathcal{G}$ is isomorphic to $\mathcal{H}_s^t(n,a,b)$.
\end{corollary}

\subsection*{Organization}

    The rest of this paper is organized as follows. Section \ref{section: pre} contains preliminary definitions and key lemmas. The proofs of our main results, Theorems \ref{thm: (1,t) CI}, \ref{thm: (s,t) CI}, \ref{thm: initial CI} and \ref{thm: initial st CI} are presented in Sections \ref{section: Proof of (1,t)-Theorem}, \ref{section: Proof of (s,t)-Theorem}, \ref{section: initial CI} and \ref{section: initial st CI} respectively. 

\section{Notations and preliminaries}
\label{section: pre}

    Let $A,B\subset [n]$ be two sets with $A \cap B =\emptyset$. Set
    \[
        \mathcal{F}(A) 
        := \left\{F \backslash A \in \binom{[n]}{a-|A|}: A \subset F\in \mathcal{F} \right\},
    \]
    \[
        \mathcal{F}(\overline{A}) 
        := \left\{F \in \binom{[n]}{a}: F \in \mathcal{F}, F \cap A = \emptyset \right\},
    \]
    \[
        \mathcal{F}(A, \overline{B}) 
        := \left\{ F \backslash A\in \binom{[n]}{a-|A|}: A \subset F\in \mathcal{F}, F \cap B = \emptyset \right\}. 
    \]
    For $\mathcal{F} \subset \binom{[n]}{a}$ and $T \subset [n]$, set
    \[
        \widetilde{\mathcal{F}}(T) := \{ F \backslash T: F \in \mathcal{F} \}.
    \]
    Here we do not require $T\subset F$. By $\mathcal{F}, \mathcal{G}$ being cross-intersecting, we have $\tau (\mathcal{F}) \leq b$ and $\tau (\mathcal{G}) \leq a$.

    For a family $\mathcal{F} \subset 2^{[n]}$, define the {\em covers} (or {\em transversals}) of $\mathcal{F}$ by 
    \[
    \mathcal{T}(\mathcal{F}) := \{ T \in [n]: F \cap T \ne \emptyset \text{ for all } F \in \mathcal{F}\}.
    \]
    Then $\tau(\mathcal{F}) = \min \{|T|:T \in \mathcal{T}(\mathcal{F}) \}$. Let
    \[
    \mathcal{T}^{(i)}(\mathcal{F})  := \{T \in \mathcal{T}(\mathcal{F}): |T| = i \}.
    \]

\begin{proposition} \label{prosition: Mt eq}
The size of $\mathcal{M}^t(n,a,b)$ is
\begin{align*}
    |\mathcal{M}^t(n,a,b)| &= 
    \binom{n-t}{a-t} + \sum_{\ell=1}^t\left(\binom{t}{\ell}\left(\binom{n-t}{a-t+\ell} + \sum_{i=1}^\ell(-1)^i\binom{\ell}{i}\binom{n-t-i(b-1)}{a-t+\ell}\right)\right) + t\\
    &= \binom{n}{a} + \sum_{i=1}^t(-1)^i\binom{t}{i}\binom{n-ib}{a} + t.
\end{align*}
\end{proposition}

\begin{proof}
Choose $t$ integers $b_1,...,b_t$ with $b_i \in B_i$. Denote by $B' = \{b_1,...,b_t\}.$ 
For $P \subset B'$, define 
$$
\mathcal{F}_P
:= \left\{ F \in \binom{[n] \backslash B'}{a}: F \cap B' = B' \backslash P, F \cap B_i \ne \emptyset, i \in [1,t] \right\}.
$$
Then we have
  \[
      |\mathcal{F}_P| = \binom{n-t}{a-t+|P|} + \sum_{i=1}^{|P|}(-1)^i\binom{|P|}{i}\binom{n-t-i(b-1)}{a-t+|P|}.
  \]
It is not difficult to check that $\mathcal{M}^t(n,a,b) =( \sqcup_{P \subset B'} \mathcal{F}_P) \cup \{B_1,...,B_t\}$. 
This implies that
\begin{align*}
    |\mathcal{M}^t(n,a,b)| &= \sum_{P \subset B'}|\mathcal{F}_P| + t \\
    &= \binom{n-t}{a-t} + \sum_{\ell=1}^t\left(\binom{t}{\ell}\left(\binom{n-t}{a-t+\ell} + \sum_{i=1}^\ell(-1)^i\binom{\ell}{i}\binom{n-t-i(b-1)}{a-t+\ell}\right)\right) + t \\
    &= \binom{n}{a} + \sum_{i=1}^t(-1)^i\binom{t}{i}\binom{n-ib}{a} + t,
\end{align*}
as desired.
\end{proof}

The following lemmas will be used in our proofs.

\begin{lemma}[Frankl and Wang \cite{FW23}]
    Let $n,k,i$ be positive integers. Then
    \begin{equation} \label{cnie}
        \binom{n-i}{k} \geq \left( \frac{n-k-(i-1)}{n-(i-1)} \right)^i\binom{n}{k}.
    \end{equation}
\end{lemma}

\begin{lemma}[Frankl and Wang \cite{FW25}] \label{qlecnie2}
    Let $b \geq 1$, $a \geq b + 2$ and $n \geq a+b$ be integers. Then
    \begin{equation}
        \binom{n}{a} - 2\binom{n-b}{a} + \binom{n-2b}{a} + 2 \geq \binom{n}{b},
    \end{equation}
    equality holds if and only if $n=a+b$ or $b=1,a=3$.
    When $n > a+b$, then
    \[
        \binom{n}{a} - 2\binom{n-b}{a} + \binom{n-2b}{a} \geq \binom{n}{b}.
    \]
\end{lemma}

\begin{lemma} [Lubell \cite{LYM66}]
    Let $k \geq 1, t \geq 1$ and $n \geq k+t$ be integers. For $\mathcal{F} \subset \binom{[n]}{k}$, define 
    \[
    \partial^{k+t}(\mathcal{F}) := \{F' \in \binom{[n]}{k+t}: \exists F \in \mathcal{F}, F \subset F' \}.
    \]
    Then
    \begin{equation} \label{LYMeq}
        \frac{|\mathcal{F}|}{\binom{n}{k}} \leq \frac{|\partial^{k+t}(\mathcal{F})|}{\binom{n}{k+t}}.
    \end{equation}
\end{lemma}

\begin{lemma} \label{lee}
    Suppose that $\mathcal{F} \subset \binom{[n]}{a}$, $\mathcal{G} \subset \binom{[n]}{b}$ are cross-intersecting families with $a\geq b+t$ and $\tau (\mathcal{G}) \geq t\geq 2$. For $n\geq a+b$, if $\mathcal{F} = \emptyset$, then
    \[
        |\mathcal{F}| + |\mathcal{G}| \leq |\mathcal{M}^t(n,a,b)|.
    \]
\end{lemma}

\begin{proof}
    When $n=a+b$, it is not difficult to check that $|\mathcal{F}| + |\mathcal{G}| \leq \binom{a+b}{a} = |\mathcal{M}^t(a+b,a,b)|$. When $b=1$, we have $t \leq |\mathcal{G}| \leq n$. Note that $\binom{n-t}{a-t} + t = |\mathcal{M}^t(n,a,1)| \geq n$.
    If $n > a+b$ and $b >1$, then it suffices to show 
    \begin{equation} \label{snsi}
        \binom{n}{a} + \sum_{i=1}^t(-1)^i\binom{t}{i}\binom{n-ib}{a}
        \geq \binom{n}{b}.
    \end{equation}
    We apply induction on $t$. When $t=2$, by Lemma \ref{qlecnie2} we have
    \[
        \binom{n}{a} -2\binom{n-b}{a} + \binom{n-2b}{a} \geq \binom{n}{b}.
    \]
    Assume that $t \geq 3$ and (\ref{snsi}) holds for all $2 \leq t' < t$, implying
    \[
       \binom{n}{a} + \sum_{i=1}^{t'}(-1)^i\binom{t'}{i}\binom{n-ib}{a}\geq \binom{n}{b}.
    \]
    Then we have
    \begin{align*}
        \binom{n-1}{a-1} + \sum_{i=1}^{t-1}(-1)^i\binom{t-1}{i}\binom{n-ib-1}{a-1} &\geq \binom{n-1}{b}.
    \end{align*}
    Using the formula
    \[
        \binom{n}{q} - \binom{n-p}{q} = \sum_{1 \leq i \leq p}\binom{n-i}{q-1},
    \]
    we have
    \begin{equation} \label{cbne}
        \binom{n}{a} + \sum_{i = 1}^{t} (-1)^i\binom{t}{i}\binom{n-ib}{a}
        = \sum_{1 \leq j \leq b}\left(\binom{n-j}{a-1} + \sum_{i = 1}^{t-1} (-1)^i\binom{t-1}{i}\binom{n-ib-j}{a-1} \right),
    \end{equation}
    implying
    \begin{equation}
        \binom{n}{a} + \sum_{i=1}^t(-1)^i\binom{t}{i}\binom{n-ib}{a}
        = \sum_{j=1}^b\left(\binom{n-j}{a-1} + \sum_{i=1}^{t-1}(-1)^i\binom{t-1}{i}\binom{n-ib-j}{a-1}\right).
    \end{equation}
    Now we distinguish two cases.
    \setcounter{case}{0}
    \begin{case}
        $n-2 > (a-1)+b$.
    \end{case}
    Since $n-2 > (a-1)+b, a-1 \geq b + (t-1)$, then by induction hypothesis we have
    \begin{equation}
        \binom{n}{a} + \sum_{i=1}^t(-1)^i\binom{t}{i}\binom{n-ib}{a}
        \geq \binom{n-1}{b} + \binom{n-2}{b}.
    \end{equation}
    By (\ref{cnie}) we have $\binom{n-1}{b} > (\frac{t-1}{t})\binom{n}{b}$ and $\binom{n-2}{b} \geq (\frac{(t-1)b-1}{tb-1})^2\binom{n}{b} \geq (\frac{2t-3}{2t-1})^2\binom{n}{b}$. It can be checked that $\frac{t-1}{t} + (\frac{2t-3}{2t-1})^2 \geq \frac{2}{3} + \frac{9}{25} > 1.02 > 1$ holds for $t \geq 3$. So we get 
    \[\binom{n}{a} + \sum_{i=1}^t(-1)^i\binom{t}{i}\binom{n-ib}{a} \geq \binom{n}{b},\]
    as desired.

    \begin{case}
        $n-2 = (a-1)+b$.
    \end{case}    
    Since $n-1 > (a-1)+b, a-1 \geq b+(t-1)$, then by induction hypothesis we have
    \begin{equation}
        \binom{n}{a} + \sum_{i=1}^t(-1)^i\binom{t}{i}\binom{n-ib}{a}
        \geq \binom{n-1}{b} + \binom{n-2}{b} - (t-1) + \binom{n-3}{a-1}.
    \end{equation}
    Similar to the above analysis, we have $\binom{n-1}{b} + \binom{n-2}{b} > \binom{n}{b}$. By $n - 3 = (a-1) + b-1$, we have $\binom{n-3}{a-1} = \binom{(a-1) + b-1}{a-1} > a+b-2 > t + 2b -2 > t-1$. We get
    \[
       \binom{n}{a}  + \sum_{i = 1}^{t} (-1)^i\binom{t}{i}\binom{n-ib}{a}
       \geq 
       \binom{n}{b},
    \]
    as desired.
\end{proof}

\begin{lemma} \label{icip1}
    Let $t \geq 1$, $b \geq 1$, $a \geq t$ and $n \geq tb$ be integers. Then
    \[
        |\mathcal{H}^t(n,a,b)| \leq |\mathcal{M}^t(n,a,b)|,
    \]
    and equality holds if and only if $n=a+b$ or $b=1$.
\end{lemma}

\begin{proof}
    It suffices to show
    \begin{equation} \label{nsf}
        \binom{n}{a}  + \sum_{i = 1}^{t} (-1)^i\binom{t}{i}\binom{n-ib}{a} + t
        \geq 
        \binom{n}{a} - \sum_{i = 0}^{t-1}\binom{b+t-1}{i}\binom{n-b-t+1}{a-i} + \binom{b+t-1}{b}.
    \end{equation}
    We apply induction on $n$ and $a$. For fixed $a$, the base case $n=a+b$ is not difficult to check. For fixed $n$, consider the base case $a = t$, we have
    \[
        \binom{n}{a}  + \sum_{i = 1}^{t} (-1)^i\binom{t}{i}\binom{n-ib}{a} + t = b^t + t, \text{ and}
    \]
    \[
        \binom{n}{a} - \sum_{i = 0}^{t-1}\binom{b+t-1}{i}\binom{n-b-t+1}{a-i} + \binom{b+t-1}{b} = \binom{b+t-1}{t} + \binom{b+t-1}{b} = \binom{b+t}{b}.
    \]
    For $t=1$, we have $b^t + t = \binom{b+t}{b}$ holds for any $b$.
    For $t=2$, we have $\binom{b+2}{b} = \frac{(b+2)(b+1)}{2} = \frac{b^2}{2} + \frac{3}{2}b + 1$. It is not difficult to check that $b^2 +2 \geq \frac{b^2}{2} + \frac{3}{2}b + 1$ holds for any positive integer $b$.
    
    For $t \geq 3$, we have
    \begin{align*}
        \binom{b+t}{t} 
        = \frac{(b+t)(b+t-1)...(b+1)}{t!}
        = \left(\frac{b}{t}+1\right)\cdots (b+1).
    \end{align*}
    Assume first that $b \geq 4$. Consider function $f(b) =\left(\frac{\sqrt{15}}{4}b \right)^2 -  \left(\frac{b}{2} +1 \right)(b+1)$, we have $f'(b) = \frac{7}{8}b - \frac{3}{2}$ and $f'(\frac{12}{7}) = 0$. Since $4 > \frac{12}{7}$, we have $f(b) \geq f(4) = 0$, implying $\left(\frac{b}{2} +1 \right)(b+1) \leq \left(\frac{\sqrt{15}}{4}b \right)^2$ for $b \geq 4$. Also, we have $\left(\frac{b}{s} +1 \right) < \frac{\sqrt{15}}{4}b, s \in [3,t] $. Hence
    \[
        \left(\frac{b}{t}+1\right) \left(\frac{b}{t-1}+1\right)\cdots (b+1) 
        \leq \left(\frac{\sqrt{15}}{4}b \right)^t < b^t + t.
    \]
    Then assume that $b = 3$. Consider function $g(t) = 3^t + t - \frac{(t+3)(t+2)(t+1)}{6}$. Note that 
    \[
    g'(t) = (\ln3)3^t + 1 - \frac{3t^2 + 12t +11}{6}, \quad g''(t) = (\ln3)^23^t - \frac{6t + 12}{6}, \quad g'''(t) = (\ln3)^33^t - 1. 
    \]
    When $t \geq 3$, it is not difficult to check $g'''(3) > 0$, $g''(3) > 59 - 5 > 0$ and $g'(3) > 30 - 13 >0$, implying $g(t) \geq g(3) = 30 - 20 >0$ for $t \geq 3$. 
    
    Now assume that $b=2$. Consider function $g_1(t) = 2^t + t - \frac{(t+2)(t+1)}{2}$. We have 
    \[
    g_1'(t) = (\ln 2)2^t + 1 - \frac{2t+3}{2}, \quad g_1''(t) = (\ln2)^22^t - 1.
    \]
    It follows that $g_1''(3) > 3.8-1 >0$ and $g_1'(3) > 6.5 - 4.5 >0$ implying $g_1(t) \geq g_1(3) = 11 - 10 >0$ for $t \geq 3$. 
    Now assume that $n > a+b, a > t$ and (\ref{nsf}) holds for all $a+ b \leq n' < n, t \leq a' <a$.
    
    Since 
    \[
       \binom{n}{a}  + \sum_{i = 1}^{t} (-1)^i\binom{t}{i}\binom{n-ib}{a} = \binom{n-1}{a}  + \sum_{i = 1}^{t} (-1)^i\binom{t}{i}\binom{n-ib-1}{a} + \binom{n-1}{a-1}  + \sum_{i = 1}^{t} (-1)^i\binom{t}{i}\binom{n-ib-1}{a-1},
    \]
    then by induction hypothesis we have
    \begin{align*}
        &~~~~~\binom{n}{a}  + \sum_{i = 1}^{t} (-1)^i\binom{t}{i}\binom{n-ib}{a} + t \\
        &= \binom{n-1}{a}  + \sum_{i = 1}^{t} (-1)^i\binom{t}{i}\binom{n-ib-1}{a} + t + \binom{n-1}{a-1}  + \sum_{i = 1}^{t} (-1)^i\binom{t}{i}\binom{n-ib-1}{a-1} +t -t\\
        &\geq \binom{n-1}{a} - \sum_{i = 0}^{t-1}\binom{b+t-1}{i}\binom{n-b-t}{a-i} + \binom{b+t-1}{b}
        + \binom{n-1}{a-1} - \sum_{i = 0}^{t-1}\binom{b+t-1}{i}\binom{n-b-t}{a-1-i} + \binom{b+t-1}{b} - t \\
        &= \binom{n}{a} - \sum_{i = 0}^{t-1}\binom{b+t-1}{i}\binom{n-b-t+1}{a-i} + \binom{b+t-1}{b} + \left(\binom{b+t-1}{b} - t\right),
    \end{align*}
    as desired.
\end{proof}

\begin{proposition} \label{icile}
    Let $t \geq 2$, $b \geq 1, a \geq b+t$ and $n \geq a+b$ be integers. Then
    \begin{equation} \label{pr1eq}
        \binom{n}{a} - \sum_{i=0}^{t-1} \binom{b+t-1}{i} \binom{n-b-t+1}{a-i} - \binom{n}{b} + \binom{b+t-1}{b} > 0.
    \end{equation}
\end{proposition} 

\begin{proof}
    Consider the following three families 
    \[
        \mathcal{A} := \left\{ F \in \binom{[n]}{a}: |F \cap [b+t-1]| \geq t \right\},~~
        \binom{[b+t-1]}{b},~~
        \mathcal{B} := \binom{[n]}{b}.
    \]
    Define $\mathcal{P} := \{P \in 2^{[b+t-1]}: |P| \in [0,b]$ \}. For $P \in \mathcal{P}$, denote by 
    \[
    \mathcal{A}_P := \{ F \in \mathcal{A}: P \subset F \cap [b+t-1], |F \cap [b+t-1]| = |P|+t \}, \quad \mathcal{B}_P := \{ G \in \mathcal{B}: G \cap [b+t-1] = P \}
    \]
    
    One can see that $\mathcal{A}_{P_1} \cap \mathcal{A}_{P_2} = \emptyset, \mathcal{B}_{P_1} \cap \mathcal{B}_{P_2} = \emptyset$ for $P_1 \ne P_2$ and $\sqcup_{P \in \mathcal{P}} \mathcal{A}_P \subset \mathcal{A}$, $\sqcup_{P \in \mathcal{P}} \mathcal{B}_P = \mathcal{B}$.

    To prove $|\mathcal{A}| + \binom{b+t-1}{b} > |\mathcal{B}|$, it suffices to prove $|\mathcal{A}_P| \geq |\mathcal{B}_P|$ for all $P \in  \mathcal{P}$. 
    Assume first that $a = b+t$. Then one can see $|\mathcal{A}_P| > |\mathcal{B}_P|$.
    So assume that $a > b+t$. Let $\mathcal{A}'_P := \{ F \backslash [b+t-1]: F \in \mathcal{A}_P \}$ and $\mathcal{B}'_P := \{G \backslash [b+t-1]: B \in \mathcal{B}_P \}.$ It suffices to prove $|\mathcal{A}'_P| \geq |\mathcal{B}'_P|$. When $|P|=b$, we already have $|\mathcal{B}_P| = \binom{n-b-t+1}{b}$, so we just consider $|P| = i, i \in [0,b-1]$. For any $G' \in \mathcal{B}'_P$ we can find a $F' \in \mathcal{A}'_P$ such that $G' \subset F'$. By (\ref{LYMeq}) we have
    \[
        |\mathcal{B}'_P|\binom{n-b-t+1}{a-(t+i)-(b-i)} \leq |\mathcal{A}'_P|\binom{a-(t+i)}{b-i} = |\mathcal{A}'_P|\binom{a-(t+i)}{a-(t+i)-(b-i)}.
    \]
    Since $n \geq a+b$, we have $(n-b-t+1) - (a-t-i) = n - a-b+1+i > 0$, implying $\binom{n-b-t+1}{a-b-t} > \binom{a-t-i}{a-b-t}$. Thus we get $|\mathcal{A}_P| > |\mathcal{B}_P|$ for all $|P| \in [0,b-1]$.
\end{proof}

\section{Proof of Theorem \ref{thm: (1,t) CI} }
\label{section: Proof of (1,t)-Theorem}

If $n = a+b$, then since $G \notin \overline{\mathcal{F}}$ we have $|\mathcal{F}| + |\mathcal{G}| \leq \binom{a+b}{a} =  |\mathcal{M}^t(a+b,a,b)|$. So we assume that $n>a+b$. We present the proof by induction on $t$. By Theorems \ref{thm: (1,1) CI} and \ref{thm: (2,2) CI} we have $|\mathcal{F}| + |\mathcal{G}| \leq |\mathcal{M}^{t}(n,a,b)|, t\in[1,2]$, and the equality holds if and only if $\mathcal{F} \cup \mathcal{G}$ is isomorphic to $\mathcal{M}^{t}(n,a,b), t \in [1,2]$. Assume that $t\geq 3$ and Theorem \ref{thm: (1,t) CI} holds for all $1 \leq t' < t$. We choose $\mathcal{F}$ and $\mathcal{G}$ such that $|\mathcal{F}| + |\mathcal{G}|$ is maximum.

\begin{claim}
    $\tau(\mathcal{G})$ = t.
\end{claim}

\begin{proof}
    We apply shifting operation on $\mathcal{F}$ and $\mathcal{G}$ until $\mathcal{F}$, $\mathcal{G}$ are initial or $\tau(\mathcal{G}) = t$. If $\mathcal{F}$ and $\mathcal{G}$ are initial with $\tau(\mathcal{G}) > t$, then by Theorem \ref{thm: initial CI} and Lemma \ref{icip1} we have $|\mathcal{F}| + |\mathcal{G}| \leq \mathcal{H}(a,b,t) < |\mathcal{M}(a,b,t)|$, a contradiction. 
\end{proof}

Denote by $T_0 := \{a_1,...,a_t\}$ a covering set of $\mathcal{G}$. Note that
\begin{equation}
    |\mathcal{F}(T_0)| + |\mathcal{G}(\overline{T_0})| \leq \binom{n-t}{a-t}.
\end{equation}

\begin{claim}
    $\mathcal{G}(a_i,\overline{T_0 \backslash a_i }) \ne \emptyset$ for all $a_i \in T_0$.
\end{claim}
\begin{proof}
    If there exists some $a_i$ in $T_0$ such that $\mathcal{G}(a_i,\overline{T_0 \backslash a_i }) = \emptyset$, then we have $T_0 \backslash \{a_i\}$ is also a covering set of $\mathcal{G}$ which means $\tau(\mathcal{G}) \leq t-1$, a contradiction.
\end{proof}

Applying Theorem \ref{thm: (1,1) CI} to $\mathcal{F}(\overline{a_1},a_2,...,a_t)$ and $\mathcal{G}(a_1,\overline{a_2,...,a_t})$, we have
\begin{equation} \label{mp1}
    |\mathcal{F}(\overline{a_1},a_2,...,a_t)| + |\mathcal{G}(a_1,\overline{a_2,...,a_t})| 
    \leq \binom{n-t}{a-t+1} - \binom{n-t-b+1}{a-t+1} + 1.
\end{equation}
Denote by $R_{\ell} := \{a_1,a_2,...,a_{\ell} \} $ and $R'_{\ell} := T_0 \backslash R_{\ell}$ for $1 \leq \ell \leq t-1$. When $n - t > a-t+\ell+b-1$, we consider the families $\mathcal{F}(\overline{R_{\ell}},R'_{\ell})$ and $\mathcal{G}(R_{\ell},\overline{R'_{\ell}})$ for $\ell \geq 2$. 

\begin{claim}
It holds that 
    \begin{equation} \label{mp3}
    |\mathcal{F}(\overline{R_\ell},R'_\ell)| + |\mathcal{G}(R_\ell,\overline{R'_\ell})| \leq \binom{n-t}{a-t+\ell} + \sum_{i=1}^{\ell}(-1)^i\binom{\ell}{i}\binom{n-t-i(b-1)}{a-t+\ell}.
\end{equation}
\end{claim}

\begin{proof}
Set 
\[
    \mathcal{G}^{\ell} := \{G' \in \binom{[n] \backslash T_0 }{b-1}:\exists~G \in \mathcal{G}(R,\overline{T_0 \backslash R}) \text{ for all } R \subset R_l \text{ and } |R| \geq 1, G \subset G'\},
\]
It follows that $\tau (\mathcal{G}^{\ell}) \geq \ell$.

If $\ell>b$, it is not difficult to check $\mathcal{G}(R_\ell,\overline{R'_\ell}) = \emptyset$ and $|\mathcal{G}^\ell| \geq |\mathcal{G}(R_\ell,\overline{R'_\ell})| + \ell$. If $\ell \leq b$, then by (\ref{LYMeq}) we have
\begin{equation} \label{SPST}
    |\mathcal{G}^\ell| \geq |\mathcal{G}(R_\ell, \overline{R'_\ell})| \frac{\binom{n-t-(b-\ell)}{b-1-(b-\ell)}}{\binom{b-1}{b-1-(b-\ell)}}.
\end{equation}
Since $n > \max\{a+b, bt\}$, then every $(b-\ell)$-subset of $[n]\backslash T_0$ is contained in a $\binom{(n-t)-(b-\ell)}{b-1-(b-\ell)} \geq n-t-b+\ell \geq b+\ell$ subset of size $b-1$, $|\mathcal{G}^\ell| \geq b+\ell$. To prove $|\mathcal{G}^\ell| \geq |\mathcal{G}(R_\ell, \overline{R'_\ell})| + \ell$, we may assume that $|\mathcal{G}(R_\ell, \overline{R'_\ell})| \geq b$. From (\ref{SPST}) we can infer that
\[
   |\mathcal{G}^\ell| \geq |\mathcal{G}(R_\ell, \overline{R'_\ell})|\frac{n-t-b+2}{b-\ell+1} \geq |\mathcal{G}(R_\ell, \overline{R'_\ell})|\frac{b+2}{b-\ell+1} \geq |\mathcal{G}(R_\ell, \overline{R'_\ell})|+|\mathcal{G}(R_\ell, \overline{R'_\ell})|\frac{\ell}{b-\ell+1}.
\]
When $|\mathcal{G}(R_\ell, \overline{R'_\ell})| \geq b$, we have $|\mathcal{G}(R_\ell, \overline{R'_\ell})|\frac{\ell}{b-\ell+1} \geq \ell$, implying $|\mathcal{G}^\ell| \geq |\mathcal{G}(R_\ell, \overline{R'_\ell})| + \ell$. 
Since $a-t+\ell \geq (b-1) + \ell $, then by induction hypothesis and use Lemma \ref{lee} we have
\begin{equation} \label{mp2}
    |\mathcal{F}(\overline{R_\ell},R'_\ell)| + |\mathcal{G}^\ell| \leq \binom{n-t}{a-t+\ell} + \sum_{i=1}^{\ell}(-1)^i\binom{\ell}{i}\binom{n-t-i(b-1)}{a-t+\ell} + \ell,
\end{equation}
implying (\ref{mp3}) holds.

Assume first that $\ell =2$. If $\mathcal{G}(R_{\ell},\overline{R'_{\ell}}) \ne \emptyset$, then by Theorem \ref{thm: (1,2) CI} we have
\begin{align*}
    |\mathcal{F}(\overline{R_2},R'_2)| + |\mathcal{G}(R_2,\overline{R'_2})| &< |\mathcal{F}(\overline{R_2},R'_2)| + |\mathcal{G}(R_2,\overline{R'_2})| + 2
    \leq (\overline{R_2},R'_2)| + |\mathcal{G}^2| \\
    &< \binom{n-t}{a-t+2} + \sum_{i=1}^{2}(-1)^i\binom{2}{i}\binom{n-t-i(b-1)}{a-t+2} + 2,
\end{align*}
implying
\[
   |\mathcal{F}(\overline{R_\ell},R'_\ell)| + |\mathcal{G}(R_\ell,\overline{R'_\ell})| 
    < \binom{n-t}{a-t+\ell} + \sum_{i=1}^{\ell}(-1)^i\binom{\ell}{i}\binom{n-t-i(b-1)}{a-t+\ell}.
\]
Then assume that $\ell > 2$. If $\mathcal{G}(R_{\ell},\overline{R'_{\ell}}) \ne \emptyset$, then $|\mathcal{G}^\ell| > |\mathcal{G}(R_{\ell},\overline{R'_{\ell}})| + \ell$, also implying
\[
   |\mathcal{F}(\overline{R_\ell},R'_\ell)| + |\mathcal{G}(R_\ell,\overline{R'_\ell})| 
    < \binom{n-t}{a-t+\ell} + \sum_{i=1}^{\ell}(-1)^i\binom{\ell}{i}\binom{n-t-i(b-1)}{a-t+\ell}.
\]

When $n-t \leq a-t+\ell +b-1$, if $\mathcal{G}(R_\ell,\overline{R'_\ell}) = \emptyset$ or $b < \ell$, then it is not difficult to see that
\[
    |\mathcal{F}(\overline{R_\ell},R'_\ell)| + |\mathcal{G}(R_\ell,\overline{R'_\ell})| 
    \leq \binom{n-t}{a-t+\ell} = \binom{n-t}{a-t+\ell} + \sum_{i=1}^{\ell}(-1)^i\binom{\ell}{i}\binom{n-t-i(b-1)}{a-t+\ell}.
\]
Assume that $\mathcal{G}(R_\ell,\overline{R'_\ell}) \ne \emptyset$. Since $\mathcal{F}(\overline{R_\ell},R'_\ell)$ and $\mathcal{G}(R_\ell,\overline{R'_\ell})$ are non-empty cross-intersecting families, then by Theorem \ref{thm: (1,1) CI} we have
\begin{align*}
    |\mathcal{F}(\overline{R_\ell},R'_\ell)| + |\mathcal{G}(R_\ell,\overline{R'_\ell})| 
    &\leq \binom{n-t}{a-t+\ell} - \binom{n-t-(b-\ell)}{a-t+\ell} + 1 \\
    &< \binom{n-t}{a-t+\ell} - \binom{a-t+\ell+1}{a-t+\ell} + 1 \\
    &\leq \binom{n-t}{a-t+\ell} - b.
\end{align*}
If $n-t = a-t+\ell+b-1$, then
\[
  \binom{n-t}{a-t+\ell} - b 
  \leq \binom{n-t}{a-t+\ell} - \ell 
  = \binom{n-t}{a-t+\ell} + \sum_{i = 1}^{\ell} (-1)^i \binom{\ell}{i}\binom{n-i(b-1)}{a-t+\ell}. 
\]
If $n-t < a-t+\ell+b-1$, then
\[
  \binom{n-t}{a-t+\ell} - b 
  < \binom{n-t}{a-t+\ell} 
  = \binom{n-t}{a-t+\ell} + \sum_{i = 1}^{\ell} (-1)^i \binom{\ell}{i}\binom{n-i(b-1)}{a-t+\ell}.
\]
It follows that (\ref{mp3}) holds. 
\end{proof}

Now consider $\mathcal{F}(\overline{T_0})$ and $\mathcal{G}(T_0)$. If $n-t \leq a + b -1$, then similar to the above analysis, we have
\begin{equation}
    |\mathcal{F}(\overline{T_0})| + |\mathcal{G}(T_0)| 
    \leq \binom{n-t}{a} - t = \binom{n-t}{a} + \sum_{i=1}^{t}(-1)^i\binom{t}{i}\binom{n-t-i(b-1)}{a}.
\end{equation}
Assume that $n-t > a + b-1$, set
\[
    \mathcal{G}' := \{G \in \binom{[n] \backslash T_0 }{b-1}:\exists H \in \mathcal{G}(R) \text{ for all } R \subset T_0 \text{ and } |R| \geq 1, H \subset G\},
\]
We can see $\tau (\mathcal{G}') \geq t$.
If $t>b$, then it is not difficult to see $\mathcal{G}(T_0) = \emptyset$ and $|\mathcal{G}'| \geq |\mathcal{G}(T_0)| + t$. If $t \leq b$, then by (\ref{LYMeq}) we have
\begin{equation} \label{SPST3}
    |\mathcal{G}'| \geq |\mathcal{G}(T_0)| \frac{\binom{n-t-(b-t)}{b-1-(b-t)}}{\binom{b-1}{b-1-(b-t)}}.
\end{equation}
Since $n > a+b$, every $(b-t)$-subset of $[n]\backslash T_0$ is contained in $\binom{(n-t)-(b-t)}{(b-1)-(b-t)} \geq n-b > a$ subset of size $b-1$, $|\mathcal{G}'| \geq a+1$. To prove $|\mathcal{G}'| \geq |\mathcal{G}(T_0)| + t$, we may assume that $|\mathcal{G}(T_0)| \geq a-t+1 \geq b$. From (\ref{SPST3}) we can infer that
\[
    |\mathcal{G}'| \geq |\mathcal{G}(T_0)|\frac{n-b-t+2}{b-t+1} \geq |\mathcal{G}(T_0)| + |\mathcal{G}(T_0)|\frac{t}{b-t+1} \geq |\mathcal{G}(T_0)| + \frac{b-1}{b-t+1}t.
\]
Since $t \geq 2$, we have $|\mathcal{G}'| \geq |\mathcal{G}(T_0)| + t$.

Now we distinguish two cases.

\setcounter{case}{0}
\begin{case}
    $\mathcal{F}(\overline{T_0}) = \emptyset$.
\end{case}

By (\ref{snsi}), since $a \geq (b-1) + t$, we get
\begin{align*}
    |\mathcal{F}(\overline{T_0})| + |\mathcal{G}'|
    \leq \binom{n-t}{b-1} 
    < \binom{n-t}{a} + \sum_{i=1}^{t}(-1)^i\binom{t}{i}\binom{n-t-i(b-1)}{a} + t,
\end{align*}
implying
\[
    |\mathcal{F}(\overline{T_0})| + |\mathcal{G}(T_0)| 
    < \binom{n-t}{a} + \sum_{i=1}^{t}(-1)^i\binom{t}{i}\binom{n-t-i(b-1)}{a}.
\]
Then we get $|\mathcal{F}| + |\mathcal{G}| < |\mathcal{M}^t(n,a,b)|$. 

\begin{case}
    $\mathcal{F}(\overline{T_0}) \ne \emptyset.$
\end{case}

By induction on $b$. 
For the base case $b=1$, since $|\mathcal{G}| = r \geq t$, we have $|\mathcal{F}| \leq \binom{n-r}{a-r}$. 
Assume that $b\geq 2$ and $|\mathcal{F}| + |\mathcal{G}| \leq |\mathcal{M}^{t}(n,a,b)|$ holds for all $1 \leq b' < b$. 
Note that $\binom{n-r}{a-r} + r \leq \binom{n-t}{a-t} + t = |\mathcal{M}^t(n,a,1)|$.
Since $\mathcal{F}(\overline{T_0}) \subset \binom{[n]\backslash T_0}{a-t}$ and $ \mathcal{G}' \subset \binom{[n]\backslash T_0}{b-1} $ 
are cross-intersecting with $ \tau(\mathcal{G}') \geq t$ with $n-t > a + b-1$ and $ n-t \geq t(b-1)$, then by induction hypothesis we have
\begin{equation}
    |\mathcal{F}(\overline{T_0})| + |\mathcal{G}'| 
    \leq \binom{n-t}{a} + \sum_{i=1}^{t}(-1)^i\binom{t}{i}\binom{n-t-i(b-1)}{a} + t,
\end{equation}
implying
\begin{equation}
    |\mathcal{F}(\overline{T_0})| + |\mathcal{G}(T_0)| 
    \leq \binom{n-t}{a} + \sum_{i=1}^{t}(-1)^i\binom{t}{i}\binom{n-t-i(b-1)}{a}.
\end{equation}
Notice that for $t > 2$ if $\mathcal{G}(T_0) \ne \emptyset$ then $|\mathcal{G}'| > |\mathcal{G}(T_0)| + t$, implying
\[
   |\mathcal{F}(\overline{T_0})| + |\mathcal{G}(T_0)| 
    < \binom{n-t}{a} + \sum_{i=1}^{t}(-1)^i\binom{t}{i}\binom{n-t-i(b-1)}{a}.
\]

By Proposition \ref{prosition: Mt eq}, we get
\begin{align*}
    |\mathcal{F}| + |\mathcal{G}| &= \sum_{R\subset T_0}|\mathcal{F}(\overline{R}, T_0 \backslash R)| + \sum_{R\subset T_0}|\mathcal{G}(R, \overline{T_0 \backslash R})| \\
    &\leq \binom{n-t}{a-t} + \sum_{\ell=1}^t\left(\binom{t}{\ell}\left(\binom{n-t}{a-t+\ell} + \sum_{i=1}^\ell(-1)^i\binom{\ell}{i}\binom{n-t-i(b-1)}{a-t+\ell}\right)\right) + t \\
    &= \binom{n}{a} + \sum_{i=1}^t(-1)^i\binom{t}{i}\binom{n-ib}{a} + t.
\end{align*}

To achieve the maximum of $|\mathcal{F}| +|\mathcal{G}| = |\mathcal{M}^t(n,a,b)|$, we need to achieve the maximum of $|\mathcal{F}(\overline{R}, T_0 \backslash R)| + |\mathcal{G}(R, \overline{T_0 \backslash R})|$ for all $R \subset T_0$. From the discussion above, it can be concluded that 
\[
    |\mathcal{F}(\overline{R}, T_0 \backslash R)| + |\mathcal{G}(R, \overline{T_0 \backslash R})| \leq \binom{n-t}{a-t-\ell} + \sum_{i=1}^\ell(-1)^i\binom{\ell}{i}\binom{n-t-i(b-1)}{a-t+\ell} \text{ for all } R \subset T_0, |R| = \ell \geq 2,
\] 
and the inequality holds strictly when $|\mathcal{G}(R, \overline{T_0 \backslash R})| \ne \emptyset$. It can also be concluded that 
\[
    |\mathcal{F}(\overline{a_i}, T_0 \backslash a_i)| + |\mathcal{G}(a_i, \overline{T_0 \backslash a_i})| \leq \binom{n-t}{a-t+1} + \binom{n-t-b+1}{a-t+1} + 1 \text{ for all } a_i \in T_0, 
\]
when $|\mathcal{G}(a_i, \overline{T_0 \backslash a_i})| > 1$ the inequality holds strictly. Then we get $|\mathcal{F}| + |\mathcal{G}| \leq |\mathcal{M}^t(n,a,b)|$ and the equality holds if and only if $\mathcal{F}\cup \mathcal{G}$ is isomorphic to $\mathcal{M}^t(n,a,b)$. This completes the proof of Theorem \ref{thm: (1,t) CI}.

\section{Proof of Theorem \ref{thm: (s,t) CI} }
\label{section: Proof of (s,t)-Theorem}

If $n = a+b$, then since $G \notin \overline{\mathcal{F}}$ we have $|\mathcal{F}| + |\mathcal{G}| \leq \binom{a+b}{a} =  |\mathcal{M}_s^t(a+b,a,b)|$. So we assume that $n>a+b$.

Let $P := \{p_1,...,p_s \} \in \mathcal{T}^{(s)}(\mathcal{F})$. Consider the families $\mathcal{F}$ and $\mathcal{G}$ such that $|\mathcal{F}| + |\mathcal{G}|$ is maximum, implying that for all $p_i \in P$, there exists $G \in \mathcal{G}$ with $p_i \in G$. Then
\[
    |\mathcal{F}| + |\mathcal{G}| = |\cup_{i=1}^s \mathcal{F}(p_i)| + |\cup_{i=1}^s \mathcal{G}(\overline{p_i})| + |\mathcal{G}(p_1,...,p_s)| \leq |\cup_{i=1}^s \mathcal{F}(p_i)| + |\cup_{i=1}^s \mathcal{G}(\overline{p_i})| + \binom{n-s}{b-s}.
\]

\begin{claim} \label{stci c1}
It holds that
\begin{equation} \label{stci eq1}
    |\cup_{i=1}^{s} \mathcal{F}(p_i)| + |\cup_{i=1}^{s} \mathcal{G}(\overline{p_i})| \leq \sum_{j=1}^{s} \left( \binom{n-j}{a-1} + \sum_{i = 1}^{t-1} (-1)^i \binom{t-1}{i}\binom{n-ib-j}{a-1} \right) + t-1,
\end{equation}
and the inequality holds strictly if $\exists \mathcal{G}(p_i,\overline{p_j}) \ne \emptyset, i,j \in [1,s],i \ne j$ or $|\cup_{i=1}^{s} \mathcal{G}(\overline{p_i})| > t-1$.
\end{claim}

\begin{proof}
By induction on $s$.
Assume first that $s = 1$. Since $\tau(\mathcal{G}(\overline{p_i})) \geq t-1$, by Theorem \ref{thm: (1,t) CI}, we get
\[
    |\mathcal{F}(p_1)| + |\mathcal{G}(\overline{p_1})| \leq \binom{n-1}{a-1} + \sum_{i = 1}^{t-1} (-1)^i\binom{t-1}{i}\binom{n-ib-2}{a-2} + t - 1.
\]
It follows that
\[
    |\mathcal{F}| + |\mathcal{G}| = |\mathcal{F}(p_1)| + |\mathcal{G}(\overline{p_1})| + |\mathcal{G}(p_1)| \leq \binom{n-1}{a-1} + \sum_{i = 1}^{t-1} (-1)^i\binom{t-1}{i}\binom{n-ib-2}{a-2} + \binom{n-1}{b-1} + t - 1,
\]
the equality holds if and only if $\mathcal{F}\cup \mathcal{B}$ is isomorphic to $\mathcal{M}_1^t(n,a,b)$, implying 
$$
|\mathcal{F}(p_1)| + |\mathcal{G}(\overline{p_1})| \leq \binom{n-1}{a-1} + \sum_{i = 1}^{t-1} (-1)^i\binom{t-1}{i}\binom{n-ib-2}{a-2} + t - 1
$$
and the inequality holds strictly if $|\mathcal{G}(\overline{p_1})| > t-1$.

Assume that $s \geq 2$ and Claim \ref{stci c1} holds for all $1 \leq s' < s$. By induction hypothesis we have
\begin{align*}
    &~~~~|\cup_{i=1}^s \mathcal{F}(p_i)| + |\cup_{i=1}^s \mathcal{G}(\overline{p_i})| \\
    &= |\cup_{i=1}^{s-1} \mathcal{F}(p_i)| + |\cup_{i=1}^{s-1} \mathcal{G}(\overline{p_i})| + |\mathcal{F}(\overline{p_1,...,p_{s-1}},p_s)|+ |\mathcal{G}(p_1,...,p_{s-1},\overline{p_s})| \\
    &\leq \sum_{j=1}^{s-1} \left( \binom{n-j}{a-1} + \sum_{i = 1}^{t-1} (-1)^i \binom{t-1}{i}\binom{n-ib-j}{a-1} \right) + t-1 + |\mathcal{F}(\overline{p_1,...,p_{s-1}},p_s)|+ |\mathcal{G}(p_1,...,p_{s-1},\overline{p_s})|.
\end{align*} 
Now we prove 
\begin{equation} \label{stic eq1}
    |\mathcal{F}(\overline{p_1,...,p_{s-1}},p_s)|+ |\mathcal{G}(p_1,...,p_{s-1},\overline{p_s})| \leq \binom{n-s}{a-1} + \sum_{i = 1}^{t-1} (-1)^i \binom{t-1}{i}\binom{n-bi-s}{a-1}.
\end{equation}
Consider $\mathcal{F}(\overline{p_1,...,p_{s-1}},p_s) \subset \binom{[n] \backslash \{p_1,...,p_s \}}{a-1}$ and $\widetilde{\mathcal{G}(\overline{p_s})}(p_1,...,p_{s-1}) \subset \binom{[n] \backslash \{p_1,...,p_s \}}{\leq b}$, denote $\mathcal{F}(\overline{p_1,...,p_{s-1}},p_s)$ by $\mathcal{F}'$. We distinguish two cases.
\setcounter{case}{0}
\begin{case}
    $n-s > a + b-1$.
\end{case}
Set
\[
    \mathcal{G}' 
    = \left\{ G \in \binom{[n] \backslash \{p_1,...,p_s \}}{b}: \exists H \in  \widetilde{\mathcal{G}(\overline{p_s})}(p_1,...,p_{s-1}), H \subset G\right\},
\]
it is obviously that $\tau(\mathcal{G}') \geq t-1$. 

To prove (\ref{stic eq1}) it suffices to show $|\mathcal{G}'| \geq |\mathcal{G}(p_1,...,p_{s-1},\overline{p_s})| + t-1$. By (\ref{LYMeq}) we have
\begin{equation} \label{csteq1}
    |\mathcal{G}'| \geq |\mathcal{G}(p_1,...,p_{s-1},\overline{p_s})| \frac{\binom{n-s-(b-s+1)}{s-1}}{\binom{b}{s-1}}.
\end{equation}
Since $n \geq \max\{tb,a+b+1 \}$, then every $(b-s+1)$-subset of $[n]\backslash \{p_1,...,p_s\}$ is contained in $\binom{(n-s)-(b-s+1)}{b-(b-s+1)} \geq n-b-1 $ subset of size $b$. To prove $|\mathcal{G}'| \geq |\mathcal{G}(p_1,...,p_{s-1},\overline{p_s})| + t-1$, we may assume that $|\mathcal{G}(p_1,...,p_{s-1},\overline{p_s})| \geq n-b-t \geq a+b+1-b-t \geq b$. From (\ref{csteq1}) we can infer that
\begin{align*}
    |\mathcal{G}'| 
    &\geq |\mathcal{G}(p_1,...,p_{s-1},\overline{p_s})| \frac{n-b-s+1}{b-s+2} \\
    &\geq |\mathcal{G}(p_1,...,p_{s-1},\overline{p_s})| + |\mathcal{G}(p_1,...,p_{s-1},\overline{p_s})| \frac{n-2b-1}{b-s+2} \\
    &\geq |\mathcal{G}(p_1,...,p_{s-1},\overline{p_s})| + (n-2b-1).
\end{align*}
Since $n \geq a+b+1, a \geq b+t-1$, then we have $n-2b-1 \geq a+b+1-2b-1 \geq t-1$, implying $|\mathcal{G}'| \geq |\mathcal{G}(p_1,...,p_{s-1},\overline{p_s})| + t-1$.

One can see that $n-s \geq (t-1)b$ and $a-1 \geq b + (t-1)-1$. By Theorem \ref{thm: (1,t) CI} , we have
\begin{align*}
    |\mathcal{F}'| + |\mathcal{G}(p_1,...,p_{s-1},\overline{p_s})| 
    &\leq |\mathcal{F}'| + |\mathcal{G}'| - (t-1) \\
    &\leq \binom{n-s}{a-1} + \sum_{i = 1}^{t-1} (-1)^i \binom{t-1}{i}\binom{n-ib-s}{a-1} + t-1 - (t-1).
\end{align*}
If $\mathcal{G}(p_1,...,p_{s-1},\overline{p_s}) \ne \emptyset$, then by Theorem \ref{thm: (1,2) CI} we have
\begin{align*}
    |\mathcal{F}'| + |\mathcal{G}(p_1,...,p_{s-1},\overline{p_s})| 
    &\leq |\mathcal{F}'| + |\mathcal{G}'| - (t-1) \\
    &< \binom{n-s}{a-1} + \sum_{i = 1}^{t-1} (-1)^i \binom{t-1}{i}\binom{n-ib-s}{a-1} + t-1 - (t-1),
\end{align*}
the inequality holds strictly.

\begin{case}
    $n-s \leq a+b-1$.
\end{case}
If $\mathcal{G}(p_1,...,p_{s-1},\overline{p_s}) \ne \emptyset$, then, since $\mathcal{F}'$ and $\mathcal{G}(p_1,...,p_{s-1},\overline{p_s})$ are cross-intersecting families, by Theorem \ref{thm: (1,1) CI} we have
\begin{align*}
    |\mathcal{F}'| + |\mathcal{G}(p_1,...,p_{s-1},\overline{p_s})| &\leqslant
    \binom{n-s}{a-1} - \binom{n-s-(b-s+1)}{a-1} + 1 \\
    &\leq \binom{n-s}{a-1} - \binom{a}{a-1} + 1.
\end{align*}
When $n-s = a+b-1$, we have $\binom{n-s}{a-1} - (a-1) < \binom{n-s}{a-1} - (t-1) = \binom{n-s}{a-1} + \sum_{i = 1}^{t-1} (-1)^i \binom{t-1}{i}\binom{n-ib-s}{a-1}$. When $n-s < a+b-1$, we have $\binom{n-s}{a-1} - (a-1) < \binom{n-s}{a-1} = \binom{n-s}{a-1} + \sum_{i = 1}^{t-1} (-1)^i \binom{t-1}{i}\binom{n-ib-s}{a-1}$. 
It follows that 
\[
|\mathcal{F}'| + |\mathcal{G}(p_1,...,p_{s-1},\overline{p_s})| 
\leq \binom{n-s}{a-1} + \sum_{i = 1}^{t-1} (-1)^i \binom{t-1}{i}\binom{n-ib-s}{a-1},
\]
the equality holds if and only if $\mathcal{G}(p_1,...,p_{s-1},\overline{p_s}) = \emptyset$.
Now we get 
\[
    |\cup_{i=1}^s \mathcal{F}(p_i)| + |\cup_{i=1}^s \mathcal{G}(\overline{p_i})| \leq \sum_{j=1}^s \left( \binom{n-j}{a-1} + \sum_{i = 1}^{t-1} (-1)^i \binom{t-1}{i}\binom{n-ib-j}{a-1} \right) + t-1,
\]
the inequality holds strictly if there exists $\mathcal{G}(p_i,\overline{p_j}) \ne \emptyset$ for $i,j \in [1,s],i \ne j$ or $|\cup_{i=1}^s \mathcal{G}(\overline{p_i})| > t-1$. 
\end{proof}

From the analysis above, we have
\begin{align*}
    |\mathcal{F}| + |\mathcal{G}| &= |\cup_{i=1}^s \mathcal{F}(p_i)| + |\cup_{i=1}^s \mathcal{G}(\overline{p_i})| + |\mathcal{G}(p_1,...,p_s)| \\
    &\leq \sum_{j=1}^s \left( \binom{n-j}{a-1} + \sum_{i = 1}^{t-1} (-1)^i \binom{t-1}{i}\binom{n-ib-j}{a-1} \right) + \binom{n-s}{a-s} + t-1,
\end{align*}
the equality holds if and only if $\mathcal{F} \cup \mathcal{G}$ is isomorphic to $\mathcal{M}_s^t(n,a,b)$. 

When $\tau(\mathcal{F}) = s$ and $\tau(\mathcal{G})= t$, notice that we also have $|\mathcal{F}| + |\mathcal{G}| \leq |\mathcal{M}_s^t(n,a,b)|$, the equality holds if and only if $\mathcal{F} \cup \mathcal{G}$ isomorphic to $\mathcal{M}_s^t(n,a,b)$. This completes the proof of Theorem \ref{thm: (s,t) CI}.

\begin{remark}
In the case of $\tau(\mathcal{F}) = s, \tau(\mathcal{G}) \geq 1$ with $n>a+b,a \geq b$. For $a>b$, notice that $|\mathcal{F}(\overline{p_1,...,p_{\ell-1}},p_{\ell})|+ |\mathcal{G}(p_1,...,p_{\ell-1},\overline{p_{\ell}})| \leq \binom{n-\ell}{a-1}$ for all $1 \leq \ell \leq s$. Applying induction on $s$ and use the same method above we can get that $|\mathcal{F}| + |\mathcal{G}| \leq \mathcal{M}_s^1(n,a,b)$ and the equality holds if and only if $\mathcal{F}\cup \mathcal{G}$ isomorphic to $\mathcal{M}_s^1(n,a,b)$. If $a=b$, then $\mathcal{M}_s^1(n,a,b)$ may not achieve the maximum of $|\mathcal{F}| + |\mathcal{G}|$. We consider the case where $s=1$, by Theorem \ref{thm: (1,1) CI} we have $|\mathcal{F}| + |\mathcal{G}| \leq \binom{n}{a} - \binom{n-a}{a} + 1$ and $|\mathcal{M}_1^1(n,a,b)| < \binom{n}{a} - \binom{n-a}{a} + 1$.
\end{remark}

\section{Proof of Theorem \ref{thm: initial CI}}
\label{section: initial CI}

\begin{proposition}
    If $\mathcal{B} \subset \binom{[n]}{b}$ is initial and $\tau(\mathcal{B}) \geq t$, then $\binom{[b+t-1]}{b} \subset \mathcal{B}$. 
\end{proposition}
    
\begin{proof}
    By induction on $t$. When $t = 1$, it is not difficult to see $[b] \subset \mathcal{B}$. Assume that $t\geq 2$ and $\binom{[b+t'-1]}{b} \subset \mathcal{B}$ for $1 \leq t' < t$. For $t' = t-1$, by induction hypothesis we have $\binom{[b+t-2]}{b} \subset \mathcal{B}$. Since $\tau(\mathcal{B}) \geq t$, we get $\mathcal{B} \ne \binom{[b+t-2]}{b}$. Since $\mathcal{B}$ is initial, we get $\mathcal{B}' := \left\{ G \in \binom{[b+t-1]}{b}: b+t-1 \in G \right\}$ is contained in $\mathcal{B}$. Then we have $\binom{[b+t-1]}{b} \subset \mathcal{B}$.
    \end{proof}

    We present the proof by induction on $t$. When $t=1$, by Theorem \ref{thm: (1,1) CI} we have got $|\mathcal{F}| + |\mathcal{G}| \leq |\mathcal{H}^1(n,a,b)|$, the equality holds if and only if $\mathcal{F} \cup \mathcal{G}$ is isomorphic to $\mathcal{H}^1(n,a,b)$. Assume that $t \geq 2$ and $|\mathcal{F}| + |\mathcal{G}| \leq |\mathcal{H}^{t'}(n,a,b)|$, the equality holds if and only if $\mathcal{F} \cup \mathcal{G}$ is isomorphic to $\mathcal{H}^{t'}(n,a,b)$ for all $1 \leq t'<t$.
    
    Note that
    \[
         |\mathcal{F}| + |\mathcal{G}| = |\mathcal{F}(1)| + |\mathcal{G}(\overline{1})| + |\mathcal{F}(\overline{1})| + |\mathcal{G}(1)|.
    \]
    Since $\mathcal{F}$ is initial, we have $\mathcal{F}(1) \ne \emptyset$. We also have $\binom{[2,b+t-1]}{b} \subset \mathcal{G}(\overline{1})$, implying $\tau(\mathcal{G}(\overline{1})) \geq t - 1$. One can see that $\mathcal{F}(1),\mathcal{G}(\overline{1})$ are cross-intersecting families with $a-1 \geq b +(t-1)- 1$. By induction hypothesis we get
    \begin{equation} \label{icie1}
        |\mathcal{F}(1)| + |\mathcal{G}(\overline{1})| \leq 
        \binom{n-1}{a-1} - \sum_{i=0}^{t-2} \binom{b+t-2}{i} \binom{n-b-t+1}{a-1-i} + \binom{b+t-2}{b}.
    \end{equation}

    Now we distinguish two cases.
    \setcounter{case}{0}
    \begin{case}
        $\mathcal{F}(\overline{1}) \ne \emptyset$.
    \end{case}
    By induction on $b$. For the base case $b=1$, since $|\mathcal{G}| = r \geq t$, we have $|\mathcal{F}| \leq \binom{n-r}{a-r}$. Note that $\binom{n-r}{a-r} + r \leq \binom{n-t}{a-t} + t = |\mathcal{H}^t(n,a,1)|$. Assume that $b \geq 2$ and $|\mathcal{F}| + |\mathcal{G}| \leq |\mathcal{H}^t(n,a,b')|$ holds for all $1 \leq b' <b$.
    Since $|G| = b-1, G \in \mathcal{G}(1)$ and $\binom{[2,b+t-1]}{b-1} \subset \mathcal{G}(1)$, we have $\tau(\mathcal{G}(1)) \geq t$. By induction hypothesis, we have
    \begin{equation} \label{icie2}
        |\mathcal{F}(\overline{1})| + |\mathcal{G}(1)| \leqslant
        \binom{n-1}{a} - \sum_{i=0}^{t-1} \binom{b+t-2}{i} \binom{n-b-t+1}{a-i} + \binom{b+t-2}{b-1}.
    \end{equation}
    Combining (\ref{icie1}) and (\ref{icie2}) we get $|\mathcal{F}| + |\mathcal{G}| \leq |\mathcal{H}^t(n,a,b)|$.
        
    \begin{case}
        $\mathcal{F}(\overline{1}) = \emptyset$.
    \end{case} 
       
    It follows that $|\mathcal{G}(1)| \leq \binom{n-1}{b-1}$ and 
    \begin{equation} \label{icie3}
        |\mathcal{F}| + |\mathcal{G}| \leqslant
        \binom{n-1}{a-1} + \binom{n-1}{b-1} - \sum_{i=0}^{t-2} \binom{b+t-2}{i} \binom{n-b-t+1}{a-1-i} + \binom{b+t-2}{b}.
    \end{equation}
    We use $\alpha$ to denote the RHS of Inequality (\ref{icie3}). Then
    \begin{equation}
        |\mathcal{H}(a,b,t)| - \alpha = 
        \binom{n-1}{a} - \sum_{i=0}^{t-1} \binom{b+t-2}{i} \binom{n-b-t+1}{a-i} + \binom{b+t-2}{b-1} - \binom{n-1}{b-1}.
    \end{equation}
    By (\ref{pr1eq}), we get $|\mathcal{H}^t(n,a,b)| - \alpha > 0$, which completes the proof of Theorem \ref{thm: initial CI}.

\section{Proof of Theorem \ref{thm: initial st CI}}
\label{section: initial st CI}
    \begin{proposition} \label{st CI pr}
        If $\mathcal{F}$ is a star and $\tau(\mathcal{G}) \geq t$, then
        \begin{equation}
            |\mathcal{F}| + |\mathcal{G}| 
            \leq \binom{n-1}{a-1} - \sum_{i=0}^{t-2} \binom{b+t-2}{i} \binom{n-b-t+1}{a-i-1} + \binom{n-1}{b-1} + \binom{b+t-2}{b}.
        \end{equation}
    \end{proposition}

    \begin{proof}
    Notice that $|\mathcal{F}| + |\mathcal{G}| = |\mathcal{F}(1)| + |\mathcal{G}(\overline{1})| + |\mathcal{G}(1)|$ and $\tau(\mathcal{G}(\overline{1})) \geq t-1$. By Theorem \ref{thm: initial CI}, we have 
    \[
        |\mathcal{F}(1)| + |\mathcal{G}(\overline{1})| 
        \leq \binom{n-1}{a-1} - \sum_{i=0}^{t-2} \binom{b+t-2}{i} \binom{n-b-t+1}{a-i-1} + \binom{b+t-2}{b}.
    \]
    The assertion then follows from the fact that $|\mathcal{G}(1)| \leq \binom{n-1}{b-1}$.    
    \end{proof}

    In the case $n=a+b$, it not difficult to see that $|\mathcal{F}| + |\mathcal{G}| \leq \binom{a+b}{a} = |\mathcal{H}_s^t(n,a,b)|$. 
    We prove the theorem by induction on $t$ and $s$. 
    By Proposition \ref{st CI pr}, the assertion holds for $s=1$. By Theorem \ref{thm: (s,t) CI}, the assertion holds for $t=2$. 
    Assume that $s\geq 2$, $t\geq 3$ and the assertion holds for all $2 \leq t' <t $ and $1 \leq s'<s$.

    Note that
    \[
         |\mathcal{F}| + |\mathcal{G}| = |\mathcal{F}(1)| + |\mathcal{G}(\overline{1})| + |\mathcal{F}(\overline{1})| + |\mathcal{G}(1)|.
    \]
    Since $\mathcal{F}$ and $\mathcal{G}$ are initial families, we have $\binom{[2,a+s-1]}{a} \subset \mathcal{F}(\overline{1})$ and $\binom{[2,b+t-1]}{b-1} \subset \mathcal{G}(1)$, implying $\tau(\mathcal{F}(\overline{1})) = s-1$ and $\tau(\mathcal{G}(1)) \geq t $. One can see that $\mathcal{F}(\overline{1})$, $\mathcal{G}(1)$ are cross-intersecting families with $a \geq b-1 + t-1$. By induction hypothesis we get
    \begin{align*}
        ~
        & |\mathcal{F}(\overline{1})| + |\mathcal{G}(1)| \\
        & \leq \sum_{j=1}^{s-1} \left( \binom{n-j-1}{a-1} - \sum_{i=0}^{t-j-1} \binom{b+t-j-2}{i} \binom{n-b-t+1}{a-i-1} \right) + \binom{n-s}{b-s} + \binom{b+t-2}{b-1} - \binom{b+t-s-1}{b-s}.
    \end{align*}

    Consider $\mathcal{F}(1)$ and $\mathcal{G}(\overline{1})$. Since $\mathcal{G}$ is initial, we have $\binom{[2,b+t-1]}{b} \subset \mathcal{G}(\overline{1})$, implying $\tau(\mathcal{G}(\overline{1})) \geq t -1 $. One can see that $\mathcal{F}(1)$, $\mathcal{G}(\overline{1})$ are cross-intersecting families with $a -1 \geq b + (t-1) -1$. By Theorem \ref{thm: initial CI} we have
    \begin{align*}
        |\mathcal{F}(1)| + |\mathcal{G}(\overline{1})| \leq 
        \binom{n-1}{a-1} - \sum_{i=0}^{t-2} \binom{b+t-2}{i} \binom{n-b-t+1}{a-1-i} + \binom{b+t-2}{b}.
    \end{align*}

Combine the above two inequalities we get
\[
    |\mathcal{F}| + |\mathcal{G}| 
        \leq \sum_{j=1}^s \left( \binom{n-j}{a-1} - \sum_{i=0}^{t-j-1} \binom{b+t-j-1}{i} \binom{n-b-t+1}{a-i-1} \right) + \binom{n-s}{b-s} + \binom{b+t-1}{b} - \binom{b+t-s-1}{b-s},
        \]
which completes the proof of Theorem \ref{thm: initial st CI}.

\section*{Acknowledgement}

The research was supported in part by National Natural Science Foundation of China (Grant Nos. 12311540140, 12131013) and Guangdong Basic \& Applied Basic Research Foundation (Grant No. 2023A1515030208).

\end{spacing}
\end{document}